# Moment-Based Spectral Analysis of Random Graphs with Given Expected Degrees

Victor M. Preciado, *Member, IEEE,* and M. Amin Rahimian, *Student, IEEE*

**Abstract**—In this paper, we analyze the limiting spectral distribution of the adjacency matrix of a random graph ensemble, proposed by Chung and Lu, in which a given expected degree sequence $\overline{w}_n^T = (w_1^{(n)}, \ldots, w_n^{(n)})$ is prescribed on the ensemble. Let $\mathbf{a}_{i,j} = 1$ if there is an edge between the nodes $\{i, j\}$ and zero otherwise, and consider the normalized (random) adjacency matrix of the graph ensemble: $\mathbf{A}_n = \sqrt{n\rho_n}[\mathbf{a}_{i,j}]_{i,j=1}^n$, where $\rho_n = 1/\sum_{i=1}^n w_i^{(n)}$. The empirical spectral distribution of $\mathbf{A}_n$ denoted by $\mathbf{F}_n(\cdot)$ is the empirical measure putting a mass $1/n$ at each of the $n$ real eigenvalues of the symmetric matrix $\mathbf{A}_n$. Under some technical conditions on the expected degree sequence, we show that with probability one, $\mathbf{F}_n(\cdot)$ converges weakly to a deterministic distribution $F(\cdot)$. Furthermore, we fully characterize this distribution by providing explicit expressions for the moments of $F(\cdot)$. We apply our results to well known degree distributions, such as power law and exponential. The asymptotic expressions of the spectral moments in each case provide significant insights about the bulk behavior of the eigenvalue spectrum.

## 1 INTRODUCTION

UNDERSTANDING the relationship between structural and spectral properties of a network is a key question in the field of Network Science. Spectral graph methods (see [1]–[3], and references therein) have become a fundamental tool in the analysis of large complex networks, and related disciplines, with a broad range of applications in machine learning, data mining, web search and ranking, scientific computing, and computer vision. Studying the relationship between the structure of a graph and its eigenvalues is the central topic in the field of *algebraic graph theory* [1]–[5]. In particular, the eigenvalues of matrices representing the graph structure, such as the adjacency or the Laplacian matrices, have a direct connection to the behavior of several networked dynamical processes, such as spreading processes [6], synchronization of oscillators [7], random walks [8], consensus dynamics [9], and a wide variety of distributed algorithms [10].

The availability of massive databases describing a great variety of real-world networks allows researchers to explore their structural properties with great detail. Statistical analysis of empirical data has unveiled the existence of multiple common patterns in a large variety of network properties, such as power-law degree distributions [11], or the small-world phenomenon [12]. Random graphs models are the tool-of-choice to analyze the connection between structural and spectral network properties. Aiming to replicate empirical observations, a variety of synthetic network models has been proposed in the literature [11], [12]. The structural property that have (arguably) attracted the most attention is the degree distribution. Empirical studies show that the degree distribution of important real-world networks, such as the Internet [13], Facebook, or Twitter, are heavy-tailed and can be approximated using power-law distributions [14], [15].

Random graph models, such as the Erdős-Rényi [16], [17], the scale-free [15], and the small-world models [18], are a versatile tool for investigating the properties of real-world networks. We find in the literature several random graphs able to model degree distributions from empirical data. One of best known graphs is the configuration model, originally proposed by Bender and Canfield in [19]. This model is able to fit a given degree sequence *exactly* (under certain technical conditions). Although many structural properties of this model have been studied in depth [20], [21], it is not specially amenable to spectral analysis. In contrast, the preferential attachment model proposed in [15] provides a justification for the emergence of power-law degree distributions in real-world networks. A tractable alternative to the preferential attachment model was proposed by Chung and Lu in [22], and analyzed in [23]–[25]. In this model, which we refer to as the Chung-Lu model, an *expected* degree sequence is prescribed onto a random graph ensemble, that can be algebraically described using a (random) adjacency matrix.

Studies of the statistical properties of the eigenvalues of random graphs and networks are prevalent in many applied areas. Examples include the investigations of the spacing between nearest eigenvalues in random models [26], [27], as well as real-word networks [28]–[30]. Empirical observations highlight spectral features not observed in classical random matrix ensembles, such as a triangle-like eigenvalue distribution in power-law networks [31], [32] or an exponential decay in the tails of the eigenvalue distribution [33], [34].

### 1.1 Main Contributions

In this work we offer an exact characterization for the eigenvalue spectrum of the adjacency matrix of random graphs with a given expected degree sequence. This characteriza-tion is in terms of the moments of the eigenvalue distributions and it hinges upon application of the moments method from random matrix theory [26], [27]. Accordingly, we give closed-form expressions describing the almost sure limits of the spectral moments. We then use these moments to draw important conclusions about the graph spectrum and to bound several spectral quantities of interest. Analyzing the limiting spectral distribution is important when designing statistical tests to investigate the structure of a network from the observed interconnection data. Such scenarios arise frequently in community detection, where stochastic block models (SBM) find widespread use. SBM captures the different edge probabilities for different communities but fails to account for variations of node degrees [35]. Subsequently, the degree-corrected stochastic block model offers an improved and more realistic null (as compared to the Erdős-Rényi model) when testing different hypothesis for community detection [36]. The knowledge of spectrum is also important in system-theoretic analysis of networked systems with random interconnections: indeed, we have applied our methods to study the controllability Gramian spectra for such systems [37].

---

• The authors are with the Department of Electrical and Systems Engineering at the University of Pennsylvania, Philadelphia, PA 19104.

The remainder of this paper is organized as follows. Preliminaries on the background and motivation of our study, as well as the random graph model under consideration, are presented in Section 2. Our main results on the asymptotic spectral moments of the adjacencies of random graphs and the characterization of the limiting spectral distributions (for the normalized adjacencies) are presented in Section 3, where we also include an outline of the proofs (which are presented in detail in Section 7). In Section 4, we apply our results to a case where node degrees are obtained by random samples from the support of a preset function and show how our main results can be applied in analysis of the spectrum of large random graphs. We consider random graphs with an exponential degree distribution as a special case and derive the asymptotic expressions of its spectral moments. In Section 5 we consider the important case of power-law degree distribution, which is known to be a good descriptor for many real-world networks. The asymptotic expressions of the spectral moments of power-law networks allow us to analyze the bulk behavior of the eigenvalue spectrum in much greater details; in particular, we can quantify and characterize the similarities with and deviations from the triangular distribution that is reported in the literature [31], [32]. Section 6 concludes the paper.

## 2 BACKGROUND & MOTIVATION

### 2.1 Chung-Lu Random Graph Model

We consider the Chung-Lu random graph model introduced in [22] and analyzed in [23]–[25], in which an expected degree sequence given by the $n$ non-negative entries of the vector $\overline{w}_n^T = (w_1^{(n)}, \ldots, w_n^{(n)})$ is prescribed over the set of nodes, labeled by $[n]$, of the graph ensemble.[1] In this model, each random edge is realized independently of all other edges and in accordance with the probability measure $\mathbb{P}\{\cdot\}$, specified below. Let $\mathbb{E}\{\cdot\}$ and $\text{Var}\{\cdot\}$ be the expectation and variance operators corresponding to $\mathbb{P}\{\cdot\}$.[2] Following [38] we allow for self-loops. To each random graph, we associate a (random) adjacency matrix, which is a zero-one matrix with the $(i,j)$-th entry being one if, and only if, there is an edge connecting nodes $i$ and $j$. The number of edges incident to a vertex is the degree of that vertex, and by the volume of a graph we mean the sum of the degrees of its vertices. In this paper, our primary interest is in characterizing the asymptotic behavior of the eigenvalues of the random adjacency matrix as the graph size $n$ increases. We consider the distribution of these eigenvalues over the real line and characterize this distribution through its moments sequence. Accordingly, the probability of having an edge between nodes $i$ and $j$ is equal to $\rho_n w_i^{(n)} w_j^{(n)}$, where $\rho_n = 1/\sum_{i=1}^n w_i^{(n)}$ is the inverse expected volume. The adjacency relations in this random graph model are represented by an $n \times n$ real-valued, symmetric random matrix $\mathbf{A}_n = \sqrt{n\rho_n}[\mathbf{a}_{ij}^{(n)}]$, where $\mathbf{a}_{ij}^{(n)}$ are independent 0-1 random variables with $\mathbb{E}\{\mathbf{a}_{ij}^{(n)}\} = \rho_n w_i^{(n)} w_j^{(n)}$. In this paper, we assume the following sparsity condition on the degree sequence[3]: $\rho_n w_i^{(n)} w_j^{(n)} = o(1)$ for all $i, j$, which will be used to ensure that the distribution of the eigenvalues of the adjacency matrix $\mathbf{A}_n$ converges, almost surely, to a deterministic distribution that is uniquely characterized by its sequence of moments when $n \to \infty$. As a main result of this paper, we explicate the technical conditions under which this convergence property holds true (cf. Assumptions A1 to A3 below). Furthermore, we proffer explicit expressions for calculating this moment sequence. These expressions, in turn, allow us to upper or lower bound various spectral metrics of practical interests (cf. Section 4 and discussions therein).

Characterization of the convergence conditions depend critically on the behavior of the extreme values of the expected degree sequence as $n$ increases. To that end, we consider two sequences $\{\hat{w}_n : n \in \mathbb{N}\}$ and $\{\check{w}_n : n \in \mathbb{N}\}$ given by $\hat{w}_n = \max_{i \in [n]} w_i^{(n)}$ and $\check{w}_n = \min_{i \in [n]} w_i^{(n)}$ for all $n \in \mathbb{N}$. Another quantity of interest, whose evolution with the graph size $n$ plays an important role, is the second-order

---

1. Throughout this paper, $\mathbb{R}$ and $\mathbb{N}$ are the set of real and natural numbers, $\mathbb{N}_0 = \{0\} \cup \mathbb{N}$, $n \in \mathbb{N}$ is a parameter denoting the size of the random graph, $[n]$ denotes $\{1, 2, \ldots, n\}$, and $\text{card}(\mathcal{X})$ denotes the cardinality of set $\mathcal{X}$. The $n \times n$ identity matrix is denoted by $I_n$, random variables are printed in boldface, matrices are denoted by capital letters. Every vector is marked by a bar over its lower case letter.

2. Let $\bar{\mathbb{P}}_n\{\cdot\}$ be the probability measures on the finite product spaces corresponding to the $n(n+1)/2$ independent entries of the $n \times n$ real symmetric random adjacency matrix $\mathbf{A}_n$ of the Chung-Lu random graph ensemble and consider the following product measures: $\mathbb{P}_n = \bar{\mathbb{P}}_n \otimes \bar{\mathbb{P}}_{n-1} \otimes \ldots \bar{\mathbb{P}}_1$ for all $n \in \mathbb{N}$. The probability measure $\mathbb{P}\{\cdot\}$ is the probability measure on the infinite product space that is the unique extension (according to Kolmogorov's existence theorem [47, Section 36]) as $n \to \infty$ of the consistent probability measures $\mathbb{P}_n$, and $\mathbb{E}\{\cdot\}$ and $\text{Var}\{\cdot\}$ are the expectation and variance operators associated with $\mathbb{P}\{\cdot\}$. We say that an event occurs almost surely if it occurs with probability one.

3. Given three functions $f(\cdot)$ $g(\cdot)$, and $h(\cdot)$ we use the asymptotic notations $f(n) = O(g(n))$ and $f(n) = o(h(n))$ to signify the relations $\limsup_{n \to \infty} |f(n)/g(n)| < \infty$ and $\lim_{n \to \infty} |f(n)/h(n)| = 0$, respectively; in the latter case we also write $h(n) = \omega(f(n))$. We use $f(n) \asymp (g(n))$ to mean that $f(n) = (1 + o(1))g(n)$.

average degree $\tilde{d}_n$, defined in [23] as,

$$\tilde{d}_n = \frac{\sum_{i=1}^{n}\left(w_i^{(n)}\right)^2}{\sum_{i=1}^{n} w_i^{(n)}} = \rho_n \sum_{i=1}^{n}\left(w_i^{(n)}\right)^2 = \rho_n\|\overline{w}_n\|_2^2.$$

The spectral moments of the adjacency matrix are derived in terms of the limiting normalized power-sums of the degree sequence, defined as follows. For each $i \in \mathbb{N}$, we define the limiting normalized $k$-th power-sum of the expected degree sequence as: $\Lambda_k := \lim_{n\to\infty} 1/n \sum_{i=1}^{n}\left(n\rho_n w_i^{(n)}\right)^k$. For our main results to hold true, we need the expected degree sequence $\overline{w}_n^T$ to satisfy the following assumptions:

**A1 (Sparse & graphical):** $\rho_n \hat{w}_n^2 < 1, \forall n$, and $\rho_n \hat{w}_n^2 = o(1)$.
**A2 (All finite moments):** $\Lambda_k$ exists (and are finite) for all $k$.
**A3 (Controlled growth of degrees):** *(i)* $\hat{w}_n/\check{w}_n = O(\log n)$, *(ii)* $\rho_n \hat{w}_n^2 = \omega(1/n)$, *(iii)* $n\rho_n = o(1)$, *(iv)* $n\rho_n \hat{w}_n = O(1)$, *(v)* $\hat{w}_n = O(\log n)$, and *(vi)* $\tilde{d}_n = o(n/\log n)$.

Note that Assumption A1, in particular, implies that the edge probabilities are all less than one and, indeed, they are asymptotically vanishing (sparsity). Furthermore, under Assumption A3 *(iv)*, with $n\rho_n \hat{\omega}_n = O(1)$, the limiting averages $\Lambda_k$ are guaranteed to be $O(1)$ and Assumption A2 ensures that they indeed exist.

The Assumptions A1 to A3 are satisfied by several practical degree distributions of interest such as power law and exponential distributions, as we will discuss in the examples and case studies later in the paper. In order to illustrate Assumption A3, let us consider an Erdős-Rényi random graph with edge probability $p_n$, so that $w_i^{(n)} = \hat{w}_n = \check{w}_n = \tilde{d}_n = np_n$ and $\rho_n = 1/(n^2 p_n)$ for all $n$. Assumptions A1 and A3 *(ii)* require $p_n = \rho_n \hat{w}_n^2 = o(1)$ and $\omega(1/n)$, respectively (i.e. $p_n \to 0$ and $np_n \to \infty$), agree with the assumptions made in the most recent analyses concerning Erdős-Rényi random graphs under the $1/\sqrt{np_n}$ normalization [39, Theorem 1.3].

### 2.2 Spectrum of the Chung-Lu Random Graph

In 2003, by a series of results, Chung, Lu and Vu established important asymptotic properties of the spectra of the adjacency matrices of random graph models with given expected degree sequences [23]–[25]. A key result of theirs specifies the almost sure limit of the largest eigenvalue of the adjacency matrix, as follows [23, Theorems 2.1 and 2.2].

***Theorem 1 (Largest Eigenvalue of Random Graphs).*** *If $\tilde{d}_n > \sqrt{\hat{w}_n}\log n$, then with probability one the largest eigenvalue of the (unnormalized) adjacency matrix is $(1+o(1))\tilde{d}_n$; while if $\sqrt{\hat{w}_n} > \tilde{d}_n \log^2 n$, then the largest eigenvalue is almost surely $(1+o(1))\sqrt{\hat{w}_n}$.*

Moreover, for a random graph model whose expected degree distribution obeys a power law, the largest eigenvalue is with probability one less than $7\sqrt{\log n}\max\{\sqrt{\hat{w}_n}, \tilde{d}_n\}$. In [23], similar conditions are established for the almost sure convergence of the $k$-th largest eigenvalues towards the square root of the $k$-th largest expected degree. In [25], relevant results concerning the spectra of other matrices, such as the Laplacian, were derived. More recent results use the machinery of concentration inequalities to investigate the behavior of the graph spectra for random graphs with independent edges [40]–[42]. In [42], the author shows concentration of the spectral norm for the Laplacian and adjacency matrices around their expectations, under certain technical conditions. These results are improved by Chung and Radcliffe [40], who use a Chernoff-type inequality to approximate the eigenvalues by those of the expected matrix and bound the error with high probability. This error bound is $O(\sqrt{\hat{w}_n \log n})$ in [40] and it is later improved to $(2+o(1))\sqrt{\hat{w}_n}$ by [41]. The authors in [43] use free probability techniques to give a characterization of the empirical spectral distribution of random graphs with arbitrary expected degrees in terms of the free multiplicative convolution of their degree distribution with the semi-circular density; however, the techniques come short of a closed-form expression for the moments of the limiting spectral distribution, and no attempt is made to determine the exact conditions that the degree distribution should satisfy for the techniques to be applicable.

## 3 MAIN RESULT

In this paper, we offer a moment-based characterization of the limiting distribution of the eigenvalues of the Chung-Lu random graph model under Assumptions A1 to A3. The moment sequence provides a versatile tool in the spectral analysis of complex networks [44]–[46]. It is worth highlighting that, in the sparsity regime $\rho_n \hat{w}_n^2 = o(1)$ and under the $\sqrt{n\rho_n}$ normalization, the largest eigenvalue of the random adjacency matrix may escape to infinity as $n \to \infty$.[4] When investigating the limiting distribution of a sequence of distributions, it may be the case that some mass escapes to infinity, which can cause the limit distribution to not be a probability distribution (does not integrate to one); in such cases, the underlying sequence of distributions is not "tight" (cf. [47, Section 25]). However, as we show in this paper, the tightness property holds for the sequence of spectral distributions in the Chung-Lu random graph model when the growth of the degrees is controlled as in Assumption A3. This is because the mass that is escaping to infinity (finitely many largest eigenvalues) is asymptotically vanishing itself: as $n \to \infty$, the contribution of finitely many eigenvalues to the continuous limiting spectral distribution in the $\rho_n \hat{w}_n^2 = o(1)$ sparse regime and under the $\sqrt{n\rho_n}$ normalization is vanishingly small. By the same token, our results complement the characterization of the largest eigenvalue in [23]–[25], which focuses on the growth rate and concentration of the largest eigenvalue, as described in Subsection 2.2.

To describe our main results we need to introduce some terminology. Let $\lambda_1(\mathbf{A}_n) \leq \lambda_2(\mathbf{A}_n) \leq \ldots \leq \lambda_n(\mathbf{A}_n)$ be $n$ real-valued random variables representing the $n$ eigenvalues of the random matrix $\mathbf{A}_n$ ordered from the smallest to the largest. We define $\delta_x\{\cdot\}$ as the probability measure on $\mathbb{R}$ assigning unit mass to the point $x \in \mathbb{R}$ and zero elsewhere. We also define $\mathcal{L}_n\{\cdot\} = (1/n)\sum_{i=1}^{n}\delta_{\lambda_i(\mathbf{A}_n)}\{\cdot\}$ as the random probability measure on the real line that assigns a mass $1/n$ to each one of the $n$ eigenvalues of the random matrix $\mathbf{A}_n$. The corresponding distribution can be written as $\mathbf{F}_n(x)$

---

4. This is true, for example, for Erdős-Rényi random graphs with edge probability $p_n$ and $\sqrt{n\rho_n} = 1/\sqrt{np_n}$ normalization, as we demonstrate in Section 3.2.

$= \mathcal{L}_n\{(-\infty, x]\} = (1/n)\text{card}(\{i \in [n] : \lambda_i(\mathbf{A}_n) \leq x\})$, where $\text{card}(\cdot)$ is the cardinality function, and is referred to as the *empirical spectral distribution* (ESD) for the random matrix $\mathbf{A}_n$. For each $x \in \mathbb{R}$, $\mathbf{F}_n(x)$ is a random variable. Moreover, we define the $k$-th spectral moment of the random matrix $\mathbf{A}_n$ as the following real-valued random variable $\mathbf{m}_k^{(n)} = (1/n)\text{trace}(\mathbf{A}_n^k) = \int_{-\infty}^{+\infty} x^k d\,\mathbf{F}_n(x)$. Our main results establish the almost sure and weak convergence of the empirical spectral distribution of the (normalized) adjacency matrix $\mathbf{A}_n$ to a deterministic distribution $F(\cdot)$. We call this distribution the *limiting spectral distribution* (LSD) and we characterize it through its moments sequence, $m_k = \int_{-\infty}^{+\infty} x^k \, dF(x)$. According to our main result, the spectral moments of the limiting spectral distribution $F(\cdot)$ are specified in terms of the limiting power-sums of the normalized degree sequence $\Lambda_k$, as follows:

---

***Theorem 2 (Spectral Moments of $\mathbf{A}_n$).*** Consider a random graph with a given expected degree sequence $\overline{w}_n = (w_1^{(n)}, \ldots, w_n^{(n)})$, satisfying Assumptions A1 to A3. With probability one, $\mathbf{F}_n(\cdot)$ converges weakly to a deterministic distribution $F(\cdot)$, which is uniquely determined by its moments $m_k = \int_{-\infty}^{+\infty} x^k \, dF(x)$. The moments of $F(\cdot)$ are given by:

$$m_{2s} = \sum_{\overline{r} \in \mathcal{R}_s} \frac{2}{s+1} \binom{s+1}{r_1, \ldots, r_s} \Lambda_1^{r_1} \Lambda_2^{r_2} \ldots \Lambda_s^{r_s}, \quad (1)$$

$$m_{2s+1} = 0, \text{ for all } s \in \mathbb{N}_0,$$

with $\mathcal{R}_s = \{\overline{r} \in \mathbb{N}_0^s : \sum_{j=1}^s r_j = s+1, \sum_{j=1}^s j\, r_j = 2s\}$ and $\Lambda_k = \lim_{n \to \infty} 1/n \sum_{i=1}^n (n\rho_n w_i^{(n)})^k$.

---

### 3.1 Proof Outline

The detailed proof is presented in Section 7.2. Here, we provide an outline of the proof steps. The crux of the argument is in showing that for each $k \in \mathbb{N}$, the $k$-th spectral moments $\mathbf{m}_k^{(n)}$ converge almost surely to $m_k$; thence concluding by the method of moments that with probability one the empirical spectral distributions $\mathbf{F}_n(\cdot)$ converge weakly to $F(\cdot)$. To begin, we consider the centralized version of the adjacency $\mathbf{A}_n$, given by

$$\hat{\mathbf{A}}_n = \mathbf{A}_n - \mathbb{E}\{\mathbf{A}_n\} = \left[\hat{\mathbf{a}}_{ij}^{(n)}\right]. \quad (2)$$

Note that the entries $\hat{\mathbf{a}}_{ij}^{(n)}$ have zero mean and, asymptotically, a rank-one pattern of variances, given by:

$$\begin{aligned}
\text{Var}\left\{\hat{\mathbf{a}}_{ij}^{(n)}\right\} &= \mathbb{E}\left\{\left(\hat{\mathbf{a}}_{ij}^{(n)}\right)^2\right\} \\
&= \rho_n w_i^{(n)} w_j^{(n)} (1 - \rho_n w_i^{(n)} w_j^{(n)}) \\
&\asymp \sigma_i^{(n)} \sigma_j^{(n)},
\end{aligned} \quad (3)$$

since, by Assumption A1, $\rho_n w_i^{(n)} w_j^{(n)} = o(1)$. This rank-one pattern is a key property that allows us to calculate closed-form expressions of the limiting spectral moments in Theorem 2.

To proceed, we introduce some necessary notations. Similarly to $\mathbf{A}_n$, we consider the eigenvalues of $\hat{\mathbf{A}}_n$ ordered from the smallest to the largest as $\lambda_1(\hat{\mathbf{A}}_n) \leq \lambda_2(\hat{\mathbf{A}}_n) \leq \ldots \leq \lambda_n(\hat{\mathbf{A}}_n)$ and define the random variable $\hat{\mathbf{m}}_k^{(n)} = (1/n)\text{trace}(\hat{\mathbf{A}}_n^k)$ to be its $k$-th spectral moment. Also let $\bar{m}_k^{(n)} = \mathbb{E}\{\hat{\mathbf{m}}_k^{(n)}\}$ be the expected spectral moments for all $k, n$. The proof of our main result proceeds as follows.

Lemma 2, included in Section 7.1, ensures that under Assumptions A3 *(i)*, *(iii)*, *(v)*, and *(vi)* $\hat{\mathbf{m}}_k^{(n)} \asymp \mathbf{m}_k^{(n)}$, almost surely for each $k \in \mathbb{N}$. Therefore, the effect of centralization on the spectral moments is asymptotically vanishing, and both $\mathbf{A}_n$ and its centralized version $\hat{\mathbf{A}}_n$ have the same limiting spectral moments. Concentration results for functionals of random matrices with independent entries imply that the spectral moments concentrate around their expected values; hence, it suffices to calculate the limiting values of the expected spectral moments $(\lim_{n \to \infty} \bar{m}_k^{(n)})$. Finally, Lemma 3 (included in Section 7.2) provides the asymptotically exact expressions for the expected spectral moments under Assumptions A1 to A3, completing the proof following the method of moments. In the next section, we use the special case of Erdős-Rényi random graphs to elaborate on these steps in more details.

### 3.2 The Case of Erdős-Rényi Random Graphs

In Erdős-Rényi random graphs, denoted by $\mathbf{G}_{n,p}$, each edge is realized with a probability $p$, independently of other edges. This is a special case of the Chung-Lu random graph model when the expected degree sequence is given by $\overline{w}_n = (np, np, \ldots, np)$. Ever since its introduction in the late 1950s by Erdős and Rényi [16], [17], properties of this well-known class of random graphs have been extensively studied [48], [49]. Indeed, the seminal work of Füredi and Komlós [50] can be used to derive asymptotic properties of the spectra of Erdős-Rényi random graphs; with probability one, putting its largest eigenvalue at $(1 + o(1))np$ and upper-bounding the absolute values of the rest by $(2 + o(1))\sqrt{np(1-p)}$. More recently, Feige and Ofek [51] have shown that under mild conditions on $p$, the largest eigenvalue of the adjacency matrix is almost surely $pn + O(\sqrt{pn})$, and all other eigenvalues are almost surely $O(\sqrt{pn})$.

Fig. 1(a) depicts the histogram of eigenvalues for a particular realization of the (unnormalized) adjacency matrix of an Erdős and Rényi graph with $n = 1000$ nodes and $p = 0.01$. In particular, we can observe that the largest eigenvalue $\lambda_1(\mathbf{A}_n) \asymp np = 10$ is located away from the remaining eigenvalues, which are asymptotically located in the $[-2\sqrt{np}, +2\sqrt{np}] \approx [-6.325, +6.325]$ interval. Note that under the $\sqrt{n\rho_n} = 1/\sqrt{np}$ normalization, the largest eigenvalue grows as $\sqrt{np}$, whereas the bulk will be compactly supported in the $[-2, 2]$ interval. Let us also consider the centralized adjacency matrix $\hat{\mathbf{A}}_n$ defined in (2). We plot a typical realization of the eigenvalue histogram of the unnormalized, centralized adjacency matrix $\hat{\mathbf{A}}_n$ in Fig. 1(b). Notice that, as pointed out in Lemma 2 (proved in Section 7.1), the effect of centralizing the adjacency matrix $\mathbf{A}_n\,(\mathbf{w})$ is to cancel out the largest eigenvalue (moving it to zero), while the bulk of eigenvalues remains (almost) unperturbed. Subsequently, the effect of normalization on the limiting spectral distribution in the $\sqrt{n\rho_n}$-normalization regime is negligible. This can be also noticed in Fig. 1(c), where we plot the empirical spectral distributions $\mathbf{F}_n(x)$ and $\hat{\mathbf{F}}_n(x)$. Indeed, under the $\sqrt{np}$-normalization, the largest eigenvalue of $\mathbf{A}_n$ grows as $\sqrt{n}$, escaping almost surely to infinity; every other eigenvalue is almost surely $O(1)$, being asymptotically compactly supported.

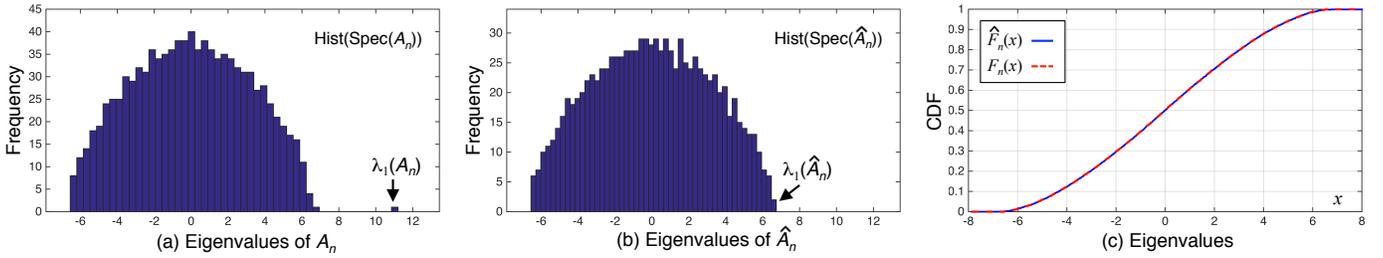

**Fig. 1:** (a) Eigenvalue histogram of the Erdős-Rényi random graph; (b) eigenvalue histogram of the centralized adjacency matrix ; (c) spectral distributions of the normalized adjacency matrix $\mathbf{A}_n$ (blue solid) and the centralized version $\hat{\mathbf{A}}_n$ (red dashed).

Indeed, as a validity check, it is possible to rederive the limiting spectral distribution of the Erdős-Rényi random graph ensemble from Theorem 2 under Assumptions A1 to A3. In particular, in this case we have that $\tilde{d}_n = \hat{w}_n = \check{w}_n = np_n$, $\rho_n = n^{-2}p_n^{-1}$ and $\Lambda_k = 1$ for all $k$; therefore, for $p_n = C \log n/n$ the random graph ensemble satisfies Assumptions A1 to A3. Theorem 2 implies that the asymptotic spectral moments satisfy

$$m_{2s} = \frac{2}{s+1} \sum_{\bar{r} \in \mathcal{R}_s} \binom{s+1}{r_1, \ldots, r_s} = \frac{1}{s+1}\binom{2s}{s} =: C_s, \quad (4)$$

These moments correspond to that of a semicircular distribution supported over $[-2, +2]$, [52], [53]. Note that to obtain a non-trivial support for the bulk of the spectrum, we need to investigate the LSD of the adjacency matrix under the normalization $\sqrt{np_n} = \sqrt{n\rho_n}$ regime. The quantity $C_s$ defined in (4) is the $s$-th Catalan number. The Catalan numbers have great significance in combinatorics: they count the total number of Dyck paths of length $2s$, cf. [26], among many other combinatorial structures. Catalan numbers can also be used to count *non-crossing partitions* of an ordered set [26], as well as many other combinatorial structures [54]. Specially relevant is the relationship between the Catalan numbers and *rooted ordered trees*. A rooted ordered tree $T$ is a tree in which one vertex is designated as the root, and the children of each vertex are ordered (see [54], page 221); i.e. there is a total order $\preceq$ on the vertex set of $T$, respecting the partial order $\preceq$ defined as follows: for all $\{j, k\} \subset \mathcal{V}(T)$, $j \preceq k$ iff $j$ belongs to the unique path on $T$ that connects $k$ to the root. There is a bijection between Dyck paths of length $2s$ and ordered trees with $s$ edges (see [26, Lemma 2.1.6]); hence, the number of ordered trees with $s$ edges is equal to $C_s$. We use $\mathcal{T}_{s+1}$ to denote the set of all rooted ordered trees on $s + 1$ vertices that are chosen without replacement from the set $[n]$. To understand the significance of the summation over the set $\mathcal{R}_s$ that appears in (4), as well as (1), we introduce some additional notation pertaing to trees and their degree sequences. Given a tree $T$ with $s$ edges and $s + 1$ vertices labeled $\{1, \ldots, s + 1\}$, the degree distribution of $T$ is defined as the sequence of integers $\bar{r}(T) = (r_1, \ldots, r_s) \in \mathbb{N}_0^s$, where $r_d = r_d(T)$ is the number of vertices with degree $d$ in $T$. We drop the tree argument $(T)$ when there is no danger of confusion: using $d_i$ for the degree of vertex $i$ and $r_d$ for number of vertices with degree $d$. Notice that the maximum degree is at most $s$; hence, $r_d = 0$ for all $d > s$. For any graph $G$ with $s$ edges, the degree distribution $(r_1, \ldots, r_s)$ satisfies $\sum_{d=1}^s d\, r_d = 2s$ and $\sum_{d=1}^s r_d = s+1$. In particular, since a connected graph with $s$ edges and $s+1$ vertices is always a tree, we have that $G$ is a tree if and only if $\sum_{d=1}^s r_d(G) = s+1$ and $\sum_{d=1}^s d\, r_d = 2s$; i.e. the set $\mathcal{R}_s$ defined in Theorem 2 denotes the set of integer sequences that are valid degree distributions for trees with $s$ edges.

We use $\mathcal{T}_{s+1}(\bar{r})$ to denote the set of all rooted ordered trees on $s + 1$ nodes whose degree distribution is $\bar{r}$. As part of the proof of the main result in (35) we show that card$(\mathcal{T}_{s+1}(\bar{r}))$, i.e. the total number of rooted ordered trees with degree distribution $\bar{r}$, is given by:

$$\text{card}\left(\mathcal{T}_{s+1}(\bar{r})\right) = \frac{2}{s+1}\binom{s+1}{r_1, \ldots, r_s}. \quad (5)$$

Next, if we use the partition $\mathcal{T}_{s+1} = \bigcup_{\bar{r} \in \mathcal{R}_s} \mathcal{T}_{s+1}(\bar{r})$, we can express the total number of ordered trees on $s+1$ vertices as, $C_s = \text{card}\left(\mathcal{T}_{s+1}\right) = \sum_{\bar{r} \in \mathcal{R}_s} \text{card}\left(\mathcal{T}_{s+1}(\bar{r})\right)$, and by replacing from (5), we obtain the right-hand side equality in (4).

## 4 Degrees Specified by Random Uniform Sampling

We now consider for Chung-Lu random graphs for which the expected degree of each node $i$ is specified as $\mathbf{w}_i^{(n)} = f_n(\mathbf{x}_i)$, where $f_n(\cdot)$ are given functions with a common support normalized to be the unit interval $[0, 1]$, and $\{\mathbf{x}_i, i \in \mathbb{N}\}$ is a random sample, uniformly and independently drawn from the unit interval.[5] To illustrate our results, let us assume that $f_n(x) = \Delta_n e^{-\alpha x}$, where $\Delta_n, \alpha_n > 0$ for all $n$, and $0 \leq x \leq 1$. Then, consider a Chung-Lu random graph with the expected degree sequence specified as $\mathbf{w}_i^{(n)} = \Delta_n e^{-\alpha \mathbf{x}_i}$, $i \in [n]$. In other words, the degrees of the resulting graph follow an exponential pattern, which have been observed in practical scenarios, such as in structural brain networks built from diffusion imaging techniques [55], [56].

The almost sure asymptotic expression for the second-order average degree $\tilde{d}_n$ can be obtained from a Monte-

---

5. Given a desired degree distribution, the function $f_n(\cdot)$ here plays the same role as the quantile function [47, Section 14] associated with the cumulative function of the desired degree distribution.

Carlo average [57, Section XVI.3], resulting in the following expression:

$$\tilde{\mathbf{d}}_n = \frac{\sum_{i=1}^n \left(\mathbf{w}_i^{(n)}\right)^2}{\sum_{i=1}^n \mathbf{w}_i^{(n)}} \asymp \frac{\int_0^1 f_n^2(x)\mathrm{d}x}{\int_0^1 f_n(x)\mathrm{d}x}$$
$$= \frac{\int_0^1 \Delta_n^2 e^{-2\alpha x}\mathrm{d}x}{\int_0^1 \Delta_n e^{-\alpha x}\mathrm{d}x} = \frac{\Delta_n(1-e^{-2\alpha})}{2(1-e^{-\alpha})}. \quad (6)$$

We know from Theorem 1 that, if $\Delta_n > \log^2 n$, then $\tilde{\mathbf{d}}_n > \sqrt{\Delta_n}\log n$, and the largest eigenvalue of the adjacency (without $\sqrt{n\rho_n}$ normalization) is almost surely given by $\tilde{\mathbf{d}}_n$. We now use Theorem 2 to write closed-form expressions for the asymptotic spectral moments of $\mathbf{A}_n$, the normalized adjacency matrix of the random graph whose expected degree sequence is given by $\overline{\mathbf{w}}_n = (\mathbf{w}_1^{(n)}, \ldots, \mathbf{w}_n^{(n)})$. We begin by calculating the almost sure asymptotic expression of the inverse expected volume of the graph, as follows,

$$\frac{1}{\boldsymbol{\rho}_n} = \sum_{i=1}^n \mathbf{w}_i^{(n)} \asymp n\int_0^1 f_n(x)\mathrm{d}x \quad (7)$$
$$= n\int_0^1 \Delta_n e^{-\alpha_n x}\mathrm{d}x = \frac{n\Delta_n}{\alpha}(1-e^{-\alpha}).$$

We now proceed to verify the qualifying conditions for applying Theorem 2. First, note the following almost sure asymptotic identities for the maximum and minimum degrees,

$$\hat{\mathbf{w}}_n = \max_{i\in[n]}\mathbf{w}_i^{(n)} \asymp \Delta_n, \quad \check{\mathbf{w}}_n = \min_{i\in[n]}\mathbf{w}_i^{(n)} \asymp \Delta_n e^{-\alpha}, \quad (8)$$

both of which can be easily verified from the corresponding order statistics for the uniform distribution on the unit interval.[6] If we set $\Delta_n = \log n$, then from the set of almost sure asymptotic identities in (6), (7) and (8) we get

$$\tilde{\mathbf{d}}_n = O(\Delta_n) = o(n/\log n), n\rho_n = O(1/\Delta_n) = o(1),$$
$$\frac{\hat{\mathbf{w}}_n}{\check{\mathbf{w}}_n} \asymp e^\alpha = O(\log n), \hat{\mathbf{w}}_n \asymp \Delta_n = O(\log n),$$
$$n\rho_n \hat{\mathbf{w}}_n = \frac{\alpha}{1-e^{-\alpha}} = O(1),$$
$$\boldsymbol{\rho}_n \hat{\mathbf{w}}_n^2 \asymp \frac{\Delta_n \alpha}{n} = \frac{\log n}{n} = o(1) \text{ and } \omega(1/n).$$

Hence, Assumptions A1 to A3 are all satisfied, and we can apply Theorem 2 to obtain the closed-form expressions for the limiting spectral moments of the normalized adjacency matrix $\mathbf{A}_n$. First note that the $k$-th order limiting averages $\boldsymbol{\Lambda}_k$ are almost surely given by

$$\boldsymbol{\Lambda}_k \asymp \int_0^1 \left(n\boldsymbol{\rho}_n f_n(x)\right)^k \mathrm{d}x = \int_0^1 \left(n\boldsymbol{\rho}_n \Delta_n e^{-\alpha_n x}\right)^k \mathrm{d}x,$$

---

6. More specifically, if we note that $\min_{i\in[n]}\mathbf{x}_i$ and $\max_{i\in[n]}\mathbf{x}_i$ are respectively distributed as $Beta(1,n)$ and $Beta(n,1)$ variables [58, Chapter 2], then the claimed almost sure limits follow by the Borel-Cantelli lemma. To see how, consider their expected values: $\mathbb{E}\{\min_{i\in[n]}\mathbf{x}_i\} = 1/(n+1)$ and $\mathbb{E}\{\max_{i\in[n]}\mathbf{x}_i\} = n/(n+1)$, and apply the Chebyshev inequality to their quadratically decaying common variance, $\mathrm{Var}\{\max_{i\in[n]}\mathbf{x}_i\} = \mathrm{Var}\{\min_{i\in[n]}\mathbf{x}_i\} = n/\left((n+1)^2(n+2)\right)$.

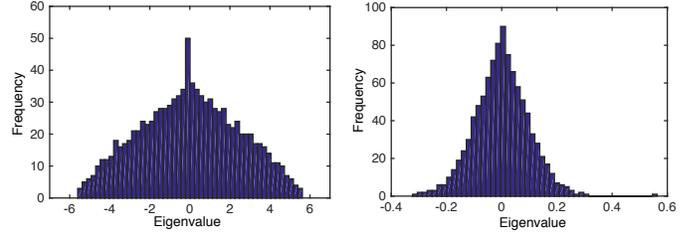

**Fig. 2:** On the left we have the eigenvalue histogram of a sample realization from the random graph model with exponential degree sequence and $n=1000$ vertices. On the right we have the eigenvalue histogram for a power-law random graph.

which results in:

$$\boldsymbol{\Lambda}_k = \frac{\alpha^{k-1}}{k(1-e^{-\alpha})^k}\left(1-e^{-k\alpha}\right).$$

Under these conditions, we can apply Theorem 2 to obtain the following closed-form expression for the asymptotic spectral moments of $\mathbf{A}_n$:

$$m_{2s} \asymp \sum_{\overline{r}\in\mathcal{R}_s} \frac{2}{s+1}\binom{s+1}{r_1,...,r_s}\prod_{k=1}^s \frac{\alpha^{(k-1)r_k}(1-e^{-k\alpha})^{r_k}}{k^{r_k}(1-e^{-\alpha})^{kr_k}}$$
$$= \frac{\alpha^{s-1}}{(1-e^{-\alpha})^{2s}}\sum_{\overline{r}\in\mathcal{R}_s}\frac{2}{s+1}\binom{s+1}{r_1,...,r_s}\prod_{k=1}^s\left(\frac{1-e^{-k\alpha}}{k}\right)^{r_k}, \quad (9)$$

where in the second equality we have used the identities $\sum_{k=1}^s r_k = s+1$, $\sum_{k=1}^s kr_k = 2s$, and $\sum_{k=1}^s(k-1)r_k = s-1$. The histogram of the eigenvalues for the normalized adjacency matrix $\mathbf{A}_n$ of a particular realization with parameters $\Delta=10$, $\alpha=1$, and $n=1000$ is plotted in Fig. 2 (left). The largest eigenvalue of the centralized adjacency matrix $\hat{\mathbf{A}}_n = \mathbf{A}_n - [\boldsymbol{\rho}_n \mathbf{w}_i^{(n)}\mathbf{w}_j^{(n)}]_{i,j=1}^n$ in this realization is given by $\lambda_n(\hat{\mathbf{A}}_n) = 5.6214$. We can upper and lower bound this eigenvalue using the $k$-th spectral moment, for any $k$, as follows [27, Equation (2.66)]:

$$\mathrm{trace}(\hat{\mathbf{A}}_n^k)^{1/k} \leq \lambda_n(\hat{\mathbf{A}}_n) \leq \left(n\cdot\mathrm{trace}(\hat{\mathbf{A}}_n^k)\right)^{1/k}. \quad (10)$$

If we consider $k=20$ in (10) and use the asymptotic spectral moment $m_{20}$ available from (9) to replace for $\mathrm{trace}(\hat{\mathbf{A}}_n^{20})$, then the lower and upper bounds on $\lambda_n(\hat{\mathbf{A}}_n)$ are given by: $(m_{20})^{1/20} = 4.3193$ and $(n\cdot m_{20})^{1/20} = 6.1011$. These values compare reasonably with the empirically observed value $\lambda_n(\hat{\mathbf{A}}_n) = 5.6214$. Furthermore, using the techniques proposed in [45], we can formulate semi-definite programs that improve the bounds in (10) by taking into account the knowledge of all spectral moments up to a fixed order, as described in [59, Section 3].

## 5 THE CASE OF POWER-LAW RANDOM GRAPHS

In this section, we study the eigenvalue distribution of the random power-law graph proposed by Chung et al. in [23]. This random graph presents an expected degree sequence given by $\overline{w}_n = (w_1^{(n)}, w_2^{(n)}, \ldots, w_n^{(n)})$ such that $w_i^{(n)} = c\, i^{-1/\beta-1}$, for $i = i_0+1,...,i_0+n$, where $\beta$ is the exponent of the power-law degree distribution; i.e. the number of nodes with degree $k$ is proportional to $k^{-\beta}$.

In this model, we can prescribe a maximum and average expected degrees, denoted by $\Delta$ and $d$, respectively, by choosing the following values of $c$ and $i_0$ [23]:

$$c = \frac{\beta-2}{\beta-1} d\, n^{\frac{1}{\beta-1}}, \quad i_0 = n \left(\frac{d(\beta-2)}{\Delta(\beta-1)}\right)^{\beta-1}.$$

For power-law degree distributions, we can asymptotically evaluate the averaged $k$-th power-sums of the expected degrees for $n \to \infty$, as follows:

$$\frac{1}{n}\sum_{i=1}^{n} w_i^k = \frac{1}{n}\sum_{i=i_0}^{i_0+n-1}\left(c\, i^{-\frac{1}{\beta-1}}\right)^k \asymp \frac{1}{n}\int_{i_0}^{i_0+n}\left(cx^{-\frac{1}{\beta-1}}\right)^k dx$$

$$= \frac{1}{n} c^k \frac{\beta-1}{\beta-1-k} x^{\frac{\beta-1-k}{\beta-1}}\Big|_{i_0}^{i_0+n} = d^k f\left(\frac{d}{\Delta}; \beta, k\right),$$

where

$$f\left(\frac{d}{\Delta}; \beta, k\right) = \left(\frac{\beta-2}{\beta-1}\right)^k \frac{\beta-1}{\beta-1-k} \times \ldots \quad (11)$$

$$\left[\left(\left(\frac{d(\beta-2)}{\Delta(\beta-1)}\right)^{\beta-1}+1\right)^{\frac{\beta-1-k}{\beta-1}} - \left(\frac{d(\beta-2)}{\Delta(\beta-1)}\right)^{\beta-1-k}\right].$$

Notice that the moments of the power-law distribution are well-defined only for $k < \beta - 1$ and the moments diverge for $k \geq \beta - 1$, [60, Section 8.4.2]. Moreover, when the largest degree is much larger than the average degree, i.e., $\Delta = \omega(d)$, the expression inside the square brackets in (11) tends to one; hence (11) simplifies as follows:

$$f\left(\frac{d}{\Delta}; \beta, k\right) = (1+o(1))\left(\frac{\beta-2}{\beta-1}\right)^k \frac{\beta-1}{\beta-1-k} = \tilde{f}(\beta, k).$$

Therefore, for $k = 1$, we have that $f\left(\frac{d}{\Delta}; \beta, k\right) \asymp 1$ and $\frac{1}{n}\sum_{i=1}^{n} w_i \asymp d$, the expected average degree. Furthermore, the above expressions can be used to compute the normalized power-sums $\Lambda_k$ in Theorem 2, which can then be used to compute closed-form values for the asymptotic expected spectral moments using (1).

***Example 1.*** *Numerical verification of the asymptotic spectral moments for the power-law degree distributions.* In the following numerical simulations, we verify the validity of Theorem 2 by computing the first five even-order spectral moments of a power-law random graph with $n = 1000$, $\beta = 3$, $\Delta = 100$, and $d = 10$. The eigenvalue histogram of one particular realization is plotted in Fig. 2 (right). In Table 1, we compare the theoretical values of the even spectral moments of the centralized adjacency matrix with the empirical values of a *single* random sample. We would like to remark how, as reported in [31], [32], the empirical spectral distribution of the power-law graph under consideration resembles a "triangular" law. In the sequel, we use Theorem 2 to study if this distribution is in fact a triangle.

## 5.1 The Triangular Spectrum

Many empirical studies of real-world networks have reported triangle-like eigenvalue spectra [31], [32]. In what follows, we want to compare the spectral density of the power-law graph with the triangular density function, given by:

$$t(x;b) = \begin{cases} \frac{2}{b} + \frac{4}{b^2}x, & \text{for } x \in [-b/2, 0], \\ \frac{2}{b} - \frac{4}{b^2}x, & \text{for } x \in (0, b/2], \\ 0, & \text{otherwise}. \end{cases} \quad (12)$$

| Order | 2 | 4 | 6 | 8 |
|---|---|---|---|---|
| Theoretical | 9.56e-3 | 3.10e-4 | 1.81e-5 | 1.46e-6 |
| Empirical | 9.5e-3 | 2.89e-4 | 1.62e-5 | 1.26e-6 |

**TABLE 1:** Theoretical versus empirical spectral moments (from one sample only) for the power-law network in Example 1.

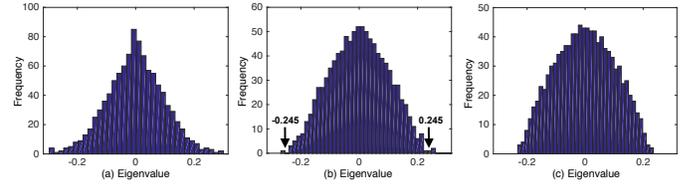

**Fig. 3:** Histograms of the eigenvalues for power-law networks with $n = 1000$, $\Delta = 100$, $d = 10$, and the following values of $\beta$: (a) $\beta = 3$, (b) $\beta = \beta_\Delta \approx 4.44$, and (c) $\beta = 6$.

The moments of this density function can be calculated as follows:

$$\widehat{m}_k(b) = \int_{-b/2}^{b/2} x^k t(x;b)\, dx \quad (13)$$

$$= \left(1 + (-1)^k\right) \frac{(b/2)^k}{(k+1)(k+2)}.$$

In what follows, we want to find conditions under which the spectral density resembles a triangular distribution. We measure the similarity in terms of the moments, in particular, in terms of the kurtosis[7] of the distribution. From (13), the second and fourth moments of the triangular law are given by $\widehat{m}_2(b) = \frac{b^2}{24}$ and $\widehat{m}_4(b) = \frac{b^4}{240}$. The kurtosis of the triangular distribution is then given by

$$\widehat{\kappa} = \frac{\widehat{m}_4(b)}{(\widehat{m}_2(b))^2} = \frac{12}{5}. \quad (14)$$

On the other hand, from Theorem 2, we can compute the kurtosis of the power-law network:

$$\kappa = \frac{m_4}{m_2^2} = \frac{2\Lambda_1^2 \Lambda_2}{(\Lambda_1^2)^2} = \frac{2\Lambda_2}{\Lambda_1^2} = 2n\rho\tilde{d}. \quad (15)$$

In our analysis, we use the difference between kurtoses as a measure of how far the spectral distribution is from the triangular law. Therefore, the spectral distribution that is closest to the triangular law is the one for which $\kappa = \widehat{\kappa}$. According to (14) and (15), $\kappa = \widehat{\kappa}$ is satisfied when $n\rho\tilde{d} = 6/5$. Furthermore, if $\kappa < 12/5$ (respectively, $\kappa > 12/5$), then the tail of the spectral distribution is "fatter" (respectively, "thiner") than the tail of the triangular law. Moreover, since $\rho \asymp \frac{1}{nd}$ and $\tilde{d} \asymp d\left(\frac{\beta-2}{\beta-1}\right)^2 \frac{\beta-1}{\beta-3}$ (when $\Delta = \omega(d)$), we have

---

7. The kurtosis of a distribution is defined as $\kappa = \mu_4/\sigma^4$, where $\mu_4$ is the fourth moment about the mean and $\sigma$ is the standard deviation. This ratio is a measurements of how heavy-tailed the distribution is.

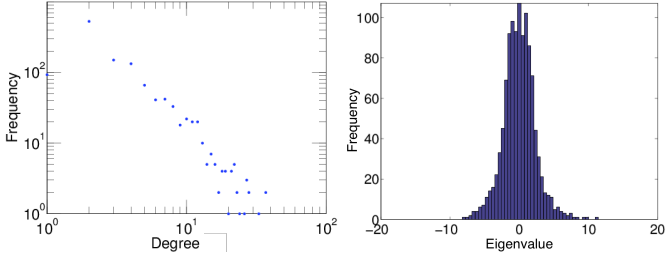

**Fig. 4:** Power-law degree distribution (left) and eigenvalue histogram (right) of the air traffic control network.

that $\kappa = \widehat{\kappa}$ when
$$\frac{(\beta-2)^2}{(\beta-1)(\beta-3)} = \frac{6}{5}.$$
Solving the above equation, we obtain the following critical value of $\beta$ for which $\kappa$ matches $\widehat{\kappa}$: $\beta_\Delta = 2 + \sqrt{6} \approx 4.44$. In other words, the kurtosis of the spectral distribution is equal to that of the triangular law only when $\beta = \beta_\Delta$.

*Remark 1.* In practice, the values of $\beta$ in real world networks are below the threshold value $\beta^* \approx 4.44$ (cf. [61]). Therefore, the spectral tails of real networks are mostly "fat", which is in accordance with empirical observations reported in [13]. For example, we include in Fig. 4 the power-law degree distribution of the (symmetric and unweighted) adjacency matrix of the air traffic control network constructed from the USA's FAA (Federal Aviation Administration) National Flight Data Center (NFDC), Preferred Routes Database [62]. Nodes in this network represent airports and links are created from strings of preferred routes recommended by the NFDC. This network has a total of 1,226 nodes, 2,615 edges, and a power-law exponent of 3.7 (which is below the threshold value $\beta_\Delta$). As predicted, the corresponding eigenvalue distribution, plotted in Fig. 4, presents a 'fat' tail.

When the spectral distribution resembles a triangular law, it is possible to approximate the support of the spectral bulk of eigenvalues by computing the value of $b$ in (12) such that the second moment of the triangular distribution matches the second moment of the theoretical spectral distribution. Since the second moment of the spectral distribution is given by $m_2 = \Lambda_1^2 = 1$ and the second moment of the triangular distribution is $\widehat{m}_2 = \frac{b^2}{24}$, the value of $b$ for which these moments match is given by: $b_\Delta = \sqrt{24}$. In Fig. 3(b), we indicate this value of $b_\Delta$ in the x-axis, which in this case is given by $n\sqrt{\rho}b_\Delta \approx 0.245$.

One may ask whether the spectral distribution of a random power-law network does indeed follow a triangular density when $\beta = \beta_\Delta$. The answer is no. This can be verified by comparing the sixth moments of the triangular distribution with $b = b_\Delta$ with the sixth moment of the power-law spectral density when $\beta = \beta_\Delta$. In particular, the sixth moment of the triangular law is given by $\widehat{m}_6(b_\Delta) = 54/7 \approx 7.71$. On the other hand, from Theorem 2, the sixth spectral moment of the power-law network is given by $m_6(\beta_\Delta) = 2\Lambda_1^3\Lambda_3 + 3\Lambda_1^2\Lambda_2^2$. Since $\Lambda_1 = 1$, $\Lambda_2 = \left(\frac{\beta_\Delta-2}{\beta_\Delta-1}\right)^2 \frac{\beta_\Delta-1}{\beta_\Delta-3} = \frac{6}{5}$, and $\Lambda_3 = \frac{6^{3/2}}{(1+\sqrt{6})^2(-2+\sqrt{6})}$, we have that
$$m_6(\beta_\Delta) = 2\left[\frac{6}{25}\left(27 + \sqrt{6}\right)\right] \approx 14.14,$$
which does not coincide with $\widehat{m}_6(b_\Delta) \approx 7.71$.

## 6 CONCLUSIONS

In this work, we have investigated the asymptotic behavior of the bulk of eigenvalues of the adjacency matrix of random graphs with given expected degree sequences. We have showed that, in the $\sqrt{n\rho_n}$ normalization regime and under some technical assumptions on the expected degrees sequence, the empirical spectral distribution of the adjacency matrix converges weakly to a deterministic distribution, which we have characterized by providing closed-form expressions for its limiting spectral moments.

We have illustrated the application of our results by analyzing the spectral distribution of large-scale networks with exponential degree distributions, which appear in structural brain networks obtained from diffusion imaging. We have also applied our results to analyze the spectrum of power-law random graphs, which are of great practical importance. Using the closed-form expressions for the asymptotic spectral moments in Theorem 2, we have investigated the triangle-like spectrum of power-law random graphs. In particular, we have provided quantitative relationships to show how the parameters of the power-law degree distribution affect the shape and properties of the graph spectrum. Furthermore, closed-form expressions of the asymptotic spectral moments allow us to bound spectral properties of practical interest, such as the support of the spectral bulk.

## 7 PROOF OF MAIN RESULTS

The argument leading to the proof of our main convergence result in Theorem 2 is based on the method of moments, and it is executed in two steps. We begin by showing (in Section 7.1, below) that with probability one $\lim_{n\to\infty} \mathbf{m}_k^{(n)} = \lim_{n\to\infty} \widehat{\mathbf{m}}_k^{(n)}$, i.e., the adjacency matrix $\mathbf{A}_n$ and its centralized version $\widehat{\mathbf{A}}_n$ share the same almost sure limits for their spectral moments. Next, in Section 7.2, we prove that these common almost sure limits are in fact given by $m_k$ in closed-form, as claimed in Theorem 2 (main result).

### 7.1 Effect of Centralization on the Spectral Moments

The following set of results measures the effect of the centralization in (2) by comparing the spectral moments of $\widehat{\mathbf{A}}_n$ and $\mathbf{A}_n$ as $n \to \infty$. Indeed, centralization by subtracting the mean $\mathbb{E}\{\mathbf{A}_n\}$ from the adjacency matrix $\mathbf{A}_n$ has the effect of shifting the largest eigenvalue towards zero. Example 1 demonstrates this shifting; however, in the sequel we shall show that the subsequent effect on the spectral moments is asymptotically vanishing, under certain mild assumptions on the degree sequences.

We know from Theorem 1 that, if $\tilde{d}_n > \sqrt{\tilde{w}_n} \log n$, then $\lambda_1(\mathbf{A}_n) \asymp \tilde{d}_n$ with probability one. We use a variation of this result to prove that the column vector of expected degrees, denoted by $\overline{w}_n^T = (w_1^{(n)}, \ldots, w_n^{(n)})^T$, is asymptotically almost surely an eigenvector of $\mathbf{A}_n$ associated with its largest eigenvalue. In particular, as $\lambda_n(\mathbf{A}_n)$ concentrates around $\tilde{d}_n$, the vector $(1/n)\mathbf{A}_n\overline{w}_n$ also concentrates around the vector $\frac{1}{n}\tilde{d}_n\overline{w}_n$. This is important when characterizing the effect of the centralization in (2) in light of the fact that

$$\mathbb{E}\{\mathbf{A}_n\} = \left[\sqrt{n}\rho_n^{3/2} w_i^{(n)} w_j^{(n)}\right]_{i,j\in[n]} = \sqrt{n}\rho_n^{3/2}\overline{w}_n\overline{w}_n^T$$
$$= \frac{\sqrt{n\rho_n}\tilde{d}_n}{\|\overline{w}_n\|_2^2}\overline{w}_n\overline{w}_n^T. \tag{16}$$

**Lemma 1 (Eigenvector Concentration).** Under Assumption A3 *(vi)*, it is true that $\frac{1}{n}\mathbf{A}_n\overline{w}_n \asymp \sqrt{\rho_n/n}\tilde{d}_n\overline{w}_n$ almost surely.

*Proof:* The $i$-th component of $\mathbf{A}_n\overline{w}_n$ is a random variable given by:

$$[\mathbf{A}_n\overline{w}_n]_i = \sum_{j=1}^n \sqrt{n\rho_n}\mathbf{a}_{ij}^{(n)} w_j^{(n)}. \tag{17}$$

Since $\mathbf{a}_{ij}^{(n)}$ is a Bernoulli random variable with $\mathbb{P}\{\mathbf{a}_{ij}^{(n)} = 1\} = \rho_n w_i^{(n)} w_j^{(n)}$, we have that

$$\mathbb{E}\{[\mathbf{A}_n\overline{w}_n]_i\} = \sqrt{n\rho_n} w_i^{(n)} \sum_{j=1}^n \rho_n \left(w_j^{(n)}\right)^2 = \sqrt{n\rho_n}\tilde{d}_n w_i^{(n)},$$

and $\mathbb{E}\{\mathbf{A}_n\overline{w}_n\} = \sqrt{n\rho_n}\tilde{d}_n\overline{w}_n$. Next, note that each of the summands $\sqrt{n\rho_n}\mathbf{a}_{ij}^{(n)} w_j^{(n)}$ in (17) are independent bounded random variables satisfying $\sqrt{n\rho_n}\mathbf{a}_{ij}^{(n)} w_j^{(n)} \in [0, \sqrt{n\rho_n} w_j^{(n)}]$ almost surely. Hence, we can apply Hoeffding's inequality [63, Theorem 2] to obtain that, for each $i$ and any $\varepsilon > 0$,

$$\mathbb{P}\left\{\frac{1}{n}\left|[\mathbf{A}_n\overline{w}_n]_i - \tilde{d}_n w_i^{(n)}\right| \geq \varepsilon\right\}$$
$$\leq 2\exp\left(\frac{-2n^2\varepsilon^2}{n\rho_n \sum_{i=1}^n \left(w_i^{(n)}\right)^2}\right) = 2e^{-2n\varepsilon^2/\tilde{d}_n}.$$

Next, note that given $\tilde{d}_n = o(n/\log n)$ per Assumption A3 *(vi)*, and for any $\alpha > 1$, we get that when $n$ is large enough $\tilde{d}_n < 2n\varepsilon^2/(\alpha\log n)$. Hence, $2e^{-2n\varepsilon^2/\tilde{d}_n} < 2/n^\alpha$ forms a summable series in $n$, and by the Borel-Cantelli lemma [47, Theorem 4.3], we get that

$$\mathbb{P}\left\{\left|\frac{1}{n}[\mathbf{A}_n\overline{w}_n]_i - \frac{\sqrt{\rho_n}\tilde{d}_n}{\sqrt{n}} w_i^{(n)}\right| \geq \varepsilon, \text{ infinitely often}\right\} = 0,$$

which holds true for any $\varepsilon > 0$. Therefore, we have

$$\mathbb{P}\left\{\lim_{n\to\infty}\frac{1}{n}[\mathbf{A}_n\overline{w}_n]_i = \lim_{n\to\infty}\frac{\sqrt{\rho_n}\tilde{d}_n}{\sqrt{n}} w_i^{(n)}\right\} = 1.$$

The claimed concentration of eigenvector around $\overline{w}_n$ now follows by the countable intersections of the above almost sure events over all $i \in \mathbb{N}$. □

We can now proceed to give conditions under which the spectral moments of $\mathbf{A}_n$ and $\hat{\mathbf{A}}_n$ are asymptotically almost surely identical, and therefore the effect of centralization on spectral moments is asymptotically vanishing.

**Lemma 2 (Vanishing Effect of Centralization).** Under Assumptions A3 *(i)*, *(iii)*, *(v)*, and *(vi)*, it is true that $\hat{\mathbf{m}}_k^{(n)} \asymp \mathbf{m}_k^{(n)}$, almost surely, for each $k \in \mathbb{N}$.

*Proof:* To begin, consider the $k=1$ case. From (2), we have

$$\hat{\mathbf{m}}_1^{(n)} = \frac{1}{n}\text{trace}(\hat{\mathbf{A}}_n) = \frac{1}{n}\text{trace}(\mathbf{A}_n - \mathbb{E}\{\mathbf{A}_n\}).$$

From (16) we know that

$$\text{trace}(\mathbb{E}\{\mathbf{A}_n\}) = \frac{\sqrt{n\rho_n}\tilde{d}_n}{\|\overline{w}_n\|_2^2}\text{trace}(\overline{w}_n\overline{w}_n^T) = \sqrt{n\rho_n}\tilde{d}_n,$$

wherefrom it follows that $\text{trace}(\hat{\mathbf{A}}_n) = \text{trace}(\mathbf{A}_n) - \sqrt{n\rho_n}\tilde{d}_n$ is true for all $n$, and in particular with probability one as $n \to \infty$. For general $k \in \mathbb{N}$, we have

$$\hat{\mathbf{m}}_k^{(n)} = \frac{1}{n}\text{trace}(\hat{\mathbf{A}}_n^k) = \frac{1}{n}\text{trace}(\mathbf{A}_n - \mathbb{E}\{\mathbf{A}_n\})^k.$$

To proceed, consider the binomial expansion of $(\mathbf{A}_n - \mathbb{E}\{\mathbf{A}_n\})^k$ consisting of a sum of the product of non-commutative elements, as follows:

$$\mathbf{A}_n^k + \mathbf{A}_n^{k-1}(-\mathbb{E}\{\mathbf{A}_n\}) + \mathbf{A}_n^{k-2}(-\mathbb{E}\{\mathbf{A}_n\})\mathbf{A}_n + \ldots$$
$$+ \mathbf{A}_n^{k-2}(-\mathbb{E}\{\mathbf{A}_n\})^2 + \mathbf{A}_n^{k-3}(-\mathbb{E}\{\mathbf{A}_n\})\mathbf{A}_n(-\mathbb{E}\{\mathbf{A}_n\}) + \ldots$$

Consider any product term of the form,

$$\mathbf{\Pi}(\overline{k}) = \mathbf{A}_n^{k_1}(-\mathbb{E}\{\mathbf{A}_n\})^{k_2}\mathbf{A}_n^{k_3}\ldots(-\mathbb{E}\{\mathbf{A}_n\})^{k_p}, \tag{18}$$

where $\overline{k} = (k_1, \ldots, k_p)$ is an integer partition of $k$, consisting of $p \geq 1$ positive integers satisfying $k_1 + \ldots + k_{p-1} = k$. Let $\tilde{k}$ be the sum of evenly indexed integers, i.e., $\tilde{k} = k_2 + k_4 + \ldots + k_{2\cdot\lfloor p/2\rfloor}$. We claim that as $n \to \infty$,

$$\frac{1}{n}\text{trace}(\mathbf{\Pi}(\overline{k})) \to \frac{(n\rho_n)^{k/2}\tilde{d}_n^k}{n}(-1)^{\tilde{k}}, \tag{19}$$

almost surely. To see why, first note that for any $k_i$, $(\overline{w}_n\overline{w}_n^T)^{k_i} = \|\overline{w}_n\|_2^{2(k_i-1)}\overline{w}_n\overline{w}_n^T$, so taking the $k_i$-th power of both sides in (16), we get

$$(-\mathbb{E}\{\mathbf{A}_n\})^{k_i} = \frac{(-\sqrt{n\rho_n}\tilde{d}_n)^{k_i}}{\|\overline{w}_n\|_2^2}\overline{w}_n\overline{w}_n^T,$$

and replacing in (18) yields

$$\mathbf{\Pi}(\overline{k}) = (-\sqrt{n\rho_n}\tilde{d}_n)^{\tilde{k}}\mathbf{A}_n^{k_1}\frac{\overline{w}_n\overline{w}_n^T}{\|\overline{w}_n\|_2^2}\mathbf{A}_n^{k_3}\ldots\frac{\overline{w}_n\overline{w}_n^T}{\|\overline{w}_n\|_2^2}. \tag{20}$$

Under Assumption A3 *(vi)*, from Lemma 1 we know that $(1/n)\mathbf{A}_n\overline{w}_n \asymp \sqrt{\rho_n/n}\tilde{d}_n\overline{w}_n$ almost surely. Multiplying both sides by $\mathbf{A}_n^{k_i-1}$ for any integer $k_i$ yields

$$\frac{1}{n}\mathbf{A}_n^{k_i}\overline{w}_n \asymp \frac{(\sqrt{n\rho_n}\tilde{d}_n)^{k_i}}{n}\overline{w}_n,$$

almost surely. Hence, taking the limits of both sides in (20) we get

$$\frac{1}{n}\mathbf{\Pi}(\overline{k}) \asymp (-1)^{\tilde{k}}\frac{(\sqrt{n\rho_n}\tilde{d}_n)^k}{n}\frac{\overline{w}_n\overline{w}_n^T}{\|\overline{w}_n\|_2^2},$$

almost surely. Taking the trace of both sides and using the fact that $\text{trace}(\overline{w}_n \overline{w}_n^T / \|\overline{w}_n\|_2^2) = 1$, leads to (19) as claimed. The above applies invariably to any of the product terms appearing in the binomial expansion of $(\mathbf{A}_n - \mathbb{E}\{\mathbf{A}_n\})^k$, with the exception of the leading term, $\text{trace}(\mathbf{A}_n^k)$, which does not simplify any further. Next, note that given any $1 \leq \tilde{k} \leq k$, the number of terms $\mathbf{\Pi}(\overline{k})$ in the binomial expansion for which (19) holds true, is exactly $\binom{k}{\tilde{k}}$; wherefore we get that with probability one,

$$\text{trace}(\hat{\mathbf{A}}_n^k) = \frac{1}{n} \text{trace}(\mathbf{A}_n - \mathbb{E}\{\mathbf{A}_n\})^k$$
$$\asymp \frac{1}{n} \left( \text{trace}(\mathbf{A}_n^k) + (\sqrt{n\rho_n}\tilde{d}_n)^k \sum_{\tilde{k}=1}^{k} (-1)^{\tilde{k}} \binom{k}{\tilde{k}} \right).$$

We next use a trite identity, $0 = (1-1)^k = 1 + \sum_{\tilde{k}=1}^{k}(-1)^{\tilde{k}} \binom{k}{\tilde{k}}$ to get

$$\frac{1}{n}\text{trace}(\hat{\mathbf{A}}_n^k) \asymp \frac{1}{n}\left(\text{trace}(\mathbf{A}_n^k) - (\sqrt{n\rho_n}\tilde{d}_n)^k\right),$$

almost surely. The claim now follows upon noting that $\frac{1}{n}(\sqrt{n\rho_n}\tilde{d}_n)^k = o(1), \forall k \in \mathbb{N}$. This is true because $\tilde{d}_n < \hat{w}_n^2/\check{w}_n$ so that $(\sqrt{n\rho_n}\tilde{d}_n)^k < (\sqrt{n\rho_n})^k (\hat{w}_n/\check{w}_n)^k \hat{w}_n^k$ and we known from Assumptions A3 *(i)* and A3 *(v)* that $(\hat{w}_n/\check{w}_n)$ and $\hat{w}_n$ are both $O(\log n)$ while Assumption A3 *(iii)* gives that $n\rho_n$ is o(1). Hence, upon division by $n$ we get that $\frac{1}{n}(\sqrt{n\rho_n}\tilde{d}_n)^k$ is $o(1)$, for all $k$ as claimed. $\square$

### 7.2 Almost Sure Limits of the Spectral Moments

Having established that the spectral moments of $\mathbf{A}_n$ and $\hat{\mathbf{A}}_n$ are asymptotically identical, in this section we establish that the spectral moments of the centralized adjacency matrix $\hat{\mathbf{A}}_n$ are indeed given by (1) in Theorem 2, i.e. we have $\hat{\mathbf{m}}_k^{(n)} \asymp m_k$, almost surely for each $k \in \mathbb{N}$. Each spectral moment is a functional that maps the $O(n^2)$ independent entries of a random matrix to a real value. Concentration of functionals of random matrices around their expected values have been established in great depth and generality. Accordingly, the deviation of these functionals (or so-called linear spectral statistics) from their expected values are characterized using Gaussian central limit theorems [64, Chapter 9], [65], as well as exponential tail bounds [26, Sections 2.3 and 4.4], [66, Section 4.14]; which, in particular, imply that $\hat{\mathbf{m}}_k^{(n)}$ converge almost surely to their limiting expected values. Hence, to determine the limiting spectral moments it suffices to calculate $\lim_{n \to \infty} \mathbb{E}\{\hat{\mathbf{m}}_k^{(n)}\} = \lim_{n \to \infty} \bar{m}_k^{(n)}$, a task which we undertake in our next lemma.

*Lemma 3 (Limiting Spectral Moments).* Under Assumptions A1 to A3, it is true that $\bar{m}_{2s}^{(n)} \asymp m_{2s}$ and $\bar{m}_{2s-1}^{(n)} = o(1)$ for all $s \in \mathbb{N}$.

*Proof:*
First, notice that

$$\bar{m}_k^{(n)} = \mathbb{E}\left\{\frac{1}{n}\text{trace}(\hat{\mathbf{A}}_n^k)\right\}$$
$$= \frac{1}{n} \sum_{1 \leq i_1,\ldots,i_k \leq n} \mathbb{E}\left\{\hat{\mathbf{a}}_{i_1 i_2}^{(n)} \hat{\mathbf{a}}_{i_2 i_3}^{(n)} \ldots \hat{\mathbf{a}}_{i_{k-1} i_k}^{(n)} \hat{\mathbf{a}}_{i_k i_1}^{(n)}\right\}, \quad (21)$$

where each summand in the last term corresponds to a sequence of $k$ nodes, which can be interpreted as a closed walk of length $k$, denoted by $w = (i_1, i_2, \ldots, i_{k-1}, i_k, i_1)$, in the complete graph with $n$ nodes, denoted by $\mathcal{K}_n$. Furthermore, we associate a weight to each walk equal to the expected value of the product of the adjacency entries along the walk. In other words, we define the weight function $\omega(w) := \mathbb{E}\{\hat{\mathbf{a}}_{i_1 i_2}^{(n)} \hat{\mathbf{a}}_{i_2 i_3}^{(n)} \ldots \hat{\mathbf{a}}_{i_{k-1} i_k}^{(n)} \hat{\mathbf{a}}_{i_k i_1}^{(n)}\}$. Our task is to study such weighted sums of closed walks. To this end, we define the set of vertices and edges visited by $w$ as $\mathcal{V}(w) = \{i_j : j \in [k]\}$ and $\mathcal{E}(w) = \{\{i_j, i_{j+1}\} : j \in [k-1]\}$. For any $e \in \mathcal{E}(w)$, we define $N(e, w)$ as the number of times that $w$ transverses the edge $e$ in any direction. We denote by $\mathcal{W}_k$ the set of all closed walks of length $k$ on $\mathcal{K}_n$ that start and end on the vertex labeled by $i_1$. It is useful to partition the set $\mathcal{W}_k$ into subsets $\mathcal{W}_{k,p}$ defined as $\mathcal{W}_{k,p} = \{w \in \mathcal{W}_k : \text{card}(\mathcal{V}(w)) = p\}$, i.e. the set of closed walks of length $k$ visiting $p$ vertices. Furthermore, it is convenient to define the following subset of $\mathcal{W}_{k,p}$:

$$\widehat{\mathcal{W}}_{k,p} = \{w \in \mathcal{W}_{k,p} : N(e, w) \geq 2 \text{ for all } e \in \mathcal{E}(w)\},$$

i.e., the set of walks in $\mathcal{W}_{k,p}$ for which each edge is traversed at least twice. Note that, any walks that do belong to $\widehat{\mathcal{W}}_{k,p}$ will have a zero weight, because for any $w \in \mathcal{W}_{k,p} \setminus \widehat{\mathcal{W}}_{k,p}$, there exists an edge $\{i, j\} \in \mathcal{E}(w)$ such that $N(\{i, j\}, w) = 1$. Hence, by independence,

$$\omega(w) = \mathbb{E}\left\{\hat{\mathbf{a}}_{ij}^{(n)}\right\} \mathbb{E}\left\{\prod_{\substack{\{k,l\} \in \\ \mathcal{E}(w) \setminus \{\{i,j\}\}}} \left(\hat{\mathbf{a}}_{kl}^{(n)}\right)^{N(\{k,l\},w)}\right\} = 0,$$

since $\mathbb{E}\{\hat{\mathbf{a}}_{ij}\} = 0$; recall from (2) that the entries $\hat{\mathbf{a}}_{ij}^{(n)}$ have zero mean. We can now employ the shorthand notation:

$$\mu_{k,p} = \frac{1}{n} \sum_{w \in \widehat{\mathcal{W}}_{k,p}} \omega(w), \quad (22)$$

to rewrite (21) as follows:

$$\bar{m}_k^{(n)} = (n\rho_n)^{k/2} \sum_{p=1}^{\lfloor k/2 \rfloor + 1} \mu_{k,p}, \quad (23)$$

where the upper limit of $p = \lfloor k/2 \rfloor + 1$ in the summation follows by the pigeonhole principle, since that every walk $w \in \mathcal{W}_{2k,p}$ with $p > k+1$ has at least one edge $e \in \mathcal{E}(w)$ such that $N(e, w) = 1$; whence, for all $p > k+1$, $\widehat{\mathcal{W}}_{2k,p} = \varnothing$.

A key step in the proof of Wigner's semicircle law [52], [53] is that for $k$ even, the summation in (23) is asymptotically dominated by those closed walks belonging to $\widehat{\mathcal{W}}_{k,k/2+1}$; i.e. those in which every edge is repeated exactly twice. In order words, for every $k$ as $n \to \infty$, $\mu_{k,k/2+1}$ dominates every other $\mu_{k,p}$ with $p < k/2 + 1$. The reason is that by repeating an edge more than twice one looses a factor of order $n$ in the available choices for vertices that comprise the closed walk, this causes the number of such walks and their contributions to be asymptotically dominated. The same principle applies to the walks that comprise the right-hand side of (23) and forms the base of our analysis.

We begin by analyzing (23) for the case $k = 2s$, i.e.

$\bar{m}_{2s}^{(n)} = (n\rho_n)^s \sum_{p=1}^{s+1} \mu_{2s,p}$. We show the desired dominance by first lower bounding the term $\mu_{2s,s+1}$ and then upper bounding the terms $\mu_{2s,p}$ for $p < s+1$, as follows. First, note that the number of walks in $\widehat{\mathcal{W}}_{2s,s+1}$ can be counted, as follows (according to Lemma 2.3.15 in [27]):

$$\text{card}(\widehat{\mathcal{W}}_{2s,s+1}) = \frac{n!}{(n-s-1)!} C_s \geq \frac{(n-s)^{s+1}}{s+1} \binom{2s}{s}. \quad (24)$$

We next lower-bound the contribution of each walk in $\widehat{\mathcal{W}}_{2s,s+1}$ as $\omega(w) \geq \rho_n^s \breve{w}_n^{2s}$, which holds true because $N(e,w) = 2, \forall e \in \mathcal{E}(w)$, and together with (24) and (22), implies that

$$\mu_{2s,s+1} \geq \frac{(n-s)^{s+1}}{n(s+1)} \binom{2s}{s} \rho_n^s \breve{w}_n^{2s}. \quad (25)$$

Note that Assumption A3 *(iii)* implies that $\rho_n = o(1/n)$ while $\breve{w}_n < \hat{w}_n = O(\log n)$; hence, $\sqrt{\rho_n} \breve{w}_n = o(\log n/\sqrt{n})$ and in particular we have that $\rho_n^s \breve{w}_n^{2s} = o(1)$ ensuring that the right hand side does not grow unbounded with increasing $n$ and is in fact vanishing with $n$. Next, to upper-bound $|\mu_{2s,p}|$ for all $p < s+1$, we make use of the following bound, which is developed by Füredi and Komlós in Section 3.2 of [50], and is subsequently used in Lemma 2 of [67] and Equation (5) of [68] as well; for all $p \in [n]$:

$$\text{card}\left(\widehat{\mathcal{W}}_{k,p}\right) \leq \frac{n!}{(n-p)!} \binom{k}{2p-2} p^{2(k-2p+2)} 2^{2p-2}$$

$$\leq n^p \binom{k}{2p-2} p^{2(k-2p+2)} 2^{2p-2}. \quad (26)$$

Furthermore, we can bound the contribution of each walk $w \in \widehat{\mathcal{W}}_{k,p}$ as $|\omega(w)| \leq \rho_n^{p-1} \hat{w}_n^{2p-2}$ because with $p$ distinct vertices in walk $w$, there is at least $p-1$ distinct edges, for each of which we can use the bound $\rho_n \hat{w}_n^2$. Indeed, for any edge $(i,j) \in \mathcal{E}(w)$, we have $N(\{i,j\},w) \geq 2$, while $|\hat{\mathbf{a}}_{ij}^{(n)}| \leq 1$ as by definition $\hat{\mathbf{a}}_{ij}^{(n)} = \mathbf{a}_{ij}^{(n)} - \rho_n w_i^{(n)} w_j^{(n)}$, so that $\hat{\mathbf{a}}_{ij}^{(n)} = 1 - \rho_n w_i^{(n)} w_j^{(n)}$ with probability $\rho_n w_i^{(n)} w_j^{(n)}$ and $\hat{\mathbf{a}}_{ij}^{(n)} = -\rho_n w_i^{(n)} w_j^{(n)}$ with probability $1 - \rho_n w_i^{(n)} w_j^{(n)}$. Using (3), we thus obtain the following bound:

$$\left|\mathbb{E}\left\{\left(\hat{\mathbf{a}}_{ij}^{(n)}\right)^{N(\{i,j\},w)}\right\}\right| \leq \mathbb{E}\left\{\left(\hat{\mathbf{a}}_{ij}^{(n)}\right)^2\right\}$$
$$= \rho_n w_i^{(n)} w_j^{(n)} (1 - \rho_n w_i^{(n)} w_j^{(n)}) \leq \rho_n \hat{w}_n^2.$$

Next, using the definition (22), together with the upper-bounds $|\omega(w)| \leq \rho_n^{p-1} \hat{w}_n^{2p-2}$ and (26), we obtain:

$$\mu_{k,p} \leq \frac{1}{n^{1-p}} \binom{k}{2p-2} p^{2(k-2p+2)} 4^{p-1} \rho_n^{p-1} \hat{w}_n^{2p-2}. \quad (27)$$

Using $k = 2s$ in (27), we get that for all $p < s+1$:

$$\mu_{2s,p} \leq \frac{1}{n^{1-p}} 4^{p-1} \binom{2s}{2p-2} p^{4(s-p+1)} \hat{w}_n^{2p-2} \rho_n^{p-1}. \quad (28)$$

To show the dominance of $\mu_{2s,s+1}$ over $\mu_{2s,p}$ for $p < s+1$, we form the ratio between the two inequalities (25) and (28) to get

$$\frac{\mu_{2s,p}}{\mu_{2s,s+1}} \leq 2\left(\frac{p^4 s^2}{4}\right)^{s+1-p} \frac{\rho_n^{p-1} \hat{w}_n^{2p-2} n^p}{\rho_n^s \breve{w}_n^{2s} (n-s)^{s+1}}$$
$$\leq \frac{s^6}{2\rho_n \hat{w}_n^2} \left(\frac{\hat{w}_n}{\breve{w}_n}\right)^{2s} \left(\frac{n}{n-s}\right)^s \frac{1}{n-s},$$

where in the first inequality we have used that

$$\binom{2s}{2p-2} \leq \left(\frac{2s}{2s-2p+2}\right)^{2s-2p+2} \leq 2s^{2s-2p+2},$$

and in the second inequality we take into account that the greatest lower bound is achieved for $p = s$, cf. [50]. The proof now follows upon noting that under Assumptions A1, A3 *(i)*, and A3 *(ii)*, we have $\rho_n \hat{w}_n^2 < 1$, $\hat{w}_n/\breve{w}_n = O(\log n)$ and $\rho_n \hat{w}_n^2 = \omega(1/n)$; whence, the rate of growth of $(s^6/\rho_n \hat{w}_n^2)(\hat{w}_n/\breve{w}_n)^{2s}$ is slower than a poly-logarithmic term times $n$, and we get that $\mu_{2s,p}/\mu_{2s,s+1} = o(1)$, as desired. Thus, having shown that $\bar{m}_{2s}^{(n)} = (n\rho_n)^s (1 + o(1)) \mu_{2s,s+1}$, we can now proceed to derive the asymptotic expressions of $\mu_{2s,s+1}$ as $n \to \infty$; in particular, to show that $\lim_{n\to\infty}(n\rho_n)^s \mu_{2s,s+1} = m_{2s}$. But before we embark on the proof of the asymptotic expressions for the even moments, we use the upper-bound in (28) to show that the odd-order moments are asymptotically vanishing: First, note from (27) that with $k = 2s+1$ each of the terms $\mu_{2s+1,p}, p \leq s+1$ can be upper-bounded, as follows:

$$|\mu_{2s+1,p}| \leq n^{p-1}\binom{2s+1}{2p-2} p^{2(2s-2p+3)} 2^{2p-2} \rho_n^{p-1} \hat{w}_n^{2p-2}$$
$$\leq \binom{2s+1}{s}(s+1)^{2(2s+3)} 2^{2s} \left(n\rho_n \hat{w}_n^2\right)^{p-1}. \quad (29)$$

Next, using (23) with $k = 2s+1$, we have $\bar{m}_{2s+1}^{(n)} = \sum_{p=1}^{s+1} \mu_{2s,p}$, and replacing from (29), the odd spectral moments can now be upper bounded as

$$|\bar{m}_{2s+1}^{(n)}|$$
$$\leq \sqrt{n\rho_n}(4n\rho_n)^s \binom{2s+1}{s}(s+1)^{4s+6} \sum_{p=1}^{s+1}\left(n\rho_n \hat{w}_n^2\right)^{p-1}$$
$$= \sqrt{n\rho_n}(4n\rho_n)^s \binom{2s+1}{s}(s+1)^{4s+6} \frac{(n\rho_n \hat{w}_n^2)^{s+1} - 1}{n\rho_n \hat{w}_n^2 - 1}$$
$$\stackrel{a}{=} \sqrt{n\rho_n} \binom{2s+1}{s}(s+1)^{4s+6} (2n\rho_n \hat{w}_n)^{2s}(1 + o(1))$$
$$\stackrel{b}{=} \sqrt{n\rho_n} \binom{2s+1}{s}(s+1)^{4s+6} O(1) \stackrel{c}{=} o(1),$$

where in ($\stackrel{a}{=}$) we use Assumption A3 *(ii)*: $\rho_n \hat{w}_n^2 = \omega(1/n)$ to conclude that $n\rho_n \hat{w}_n^2 >> 1$ and simplify the leftmost fraction; in ($\stackrel{b}{=}$) we use Assumption A3 *(iv)*: $n\rho_n \hat{w}_n = O(1)$ and in ($\stackrel{c}{=}$) we use Assumption A3 *(iii)*: $\sqrt{n\rho_n} = o(1)$, which together complete the proof of the claim about vanishing odd moments. Having thus shown the asymptotic relations $\bar{m}_{2s+1}^{(n)} = o(1)$ and $\bar{m}_{2s}^{(n)} = (1 + o(1)) \mu_{2s,s+1}$, we are now ready to show that $\mu_{2s,s+1}$ are asymptotically equal to $m_{2s}$, completing the proof of the claim about the limiting spectral moments in Lemma 3.

Starting from the definition (22) with $k = 2s$ and $p = s+1$, we have that

$$\mu_{2s,s+1} = \frac{1}{n} \sum_{w \in \widehat{\mathcal{W}}_{2s,s+1}} \omega(w) \tag{30}$$

$$= \frac{1}{n} \sum_{T \in \mathcal{T}_{s+1}} \prod_{\{i,j\} \in \mathcal{E}(T)} \rho_n w_i^{(n)} w_j^{(n)} (1 - \rho_n w_i^{(n)} w_j^{(n)}),$$

where in the last equality we write the summation in terms of the set of rooted ordered trees (using a bijection described below), and use (3) together with the independence of the entries to express the weight of each tree as the product of the terms $\rho_n w_i^{(n)} w_j^{(n)} (1 - \rho_n w_i^{(n)} w_j^{(n)})$ corresponding to each edge $\{i,j\}$ of the tree. The second equality in (30) is based on a bijection between the set of walks in $\widehat{\mathcal{W}}_{2s,s+1}$ and the set of rooted ordered trees on $s+1$ vertices, $\mathcal{T}_{s+1}$, as follows. First, note that the set of edges visited by any walk in $\widehat{\mathcal{W}}_{2s,s+1}$ is always a tree with $s+1$ vertices (since we need at least $s$ distinct edges for the induced graph to be connected). Second, it is clear that for any walk in $\widehat{\mathcal{W}}_{2s,s+1}$, each edge must be visited exactly twice and the total number of edges in the induced graph is exactly $s$. Furthermore, any walk $w \in \widehat{\mathcal{W}}_{2s,s+1}$ corresponds to a depth-first traversal of this induced tree (see [69] for more details about depth-first traversals of a tree) and we can, thus, impose a total order on the vertices of the induced tree (using the order of the first appearance in the walk $w$).

Next note that since per Assumption A1: $\rho_n w_i^{(n)} w_j^{(n)} \to 0$ as $n \to \infty$, for any $s$ fixed we obtain

$$\mu_{2s,s+1} = \frac{(1-o(1))\rho_n^s}{n} \sum_{T \in \mathcal{T}_{s+1}} \prod_{\{i,j\} \in \mathcal{E}(T)} w_i^{(n)} w_j^{(n)}. \tag{31}$$

For any tree $T$, it is always true that

$$\prod_{\{i,j\} \in \mathcal{E}(T)} w_i^{(n)} w_j^{(n)} = \prod_{i \in \mathcal{V}(T)} \left(w_i^{(n)}\right)^{d_i(T)},$$

where $d_i(T)$ is the degree of node $i$ in the tree $T$. Hence, (31) can be written as:

$$\mu_{2s,s+1} = \frac{(1-o(1))\rho_n^s}{n} \times \tag{32}$$

$$\sum_{T \in \mathcal{T}_{s+1}} \sum_{\substack{1 \le i_1, \ldots, i_{s+1} \le n \\ \text{card}(\{i_1,\ldots,i_{s+1}\}) = s+1}} \left(\prod_{k=1}^{s+1} \left(w_{i_k}^{(n)}\right)^{d_k(T)}\right).$$

Next, note that for any rooted ordered tree $T \in \mathcal{T}_{s+1}$ whose nodes are labeled by $[s+1]$ and with their respective degrees $(d_1(T), \ldots, d_{s+1}(T))$ we have:

$$\sum_{j=1}^{s+1} \sum_{\substack{1 \le i_1, \ldots, i_{s+1} \le n \\ \text{card}(\{i_1,\ldots,i_{s+1}\}) = j}} \left(\prod_{k=1}^{s+1} \left(w_{i_k}^{(n)}\right)^{d_k(T)}\right) \tag{33}$$

$$= \sum_{i_1=1}^{n} \sum_{i_2=1}^{n} \cdots \sum_{i_{s+1}=1}^{n} \left(w_{i_1}^{(n)}\right)^{d_1(T)} \cdots \left(w_{i_{s+1}}^{(n)}\right)^{d_{s+1}(T)}$$

$$= \prod_{j=1}^{s+1} \left(\sum_{k=1}^{n} \left(w_k^{(n)}\right)^{d_j(T)}\right) = \prod_{j=1}^{s+1} S_{n,d_j(T)} = \prod_{j=1}^{s} S_{n,j}^{r_j(T)}$$

where we use the partial power-sums notation $S_{n,k} = \sum_{i=1}^{n} \left(w_i^{(n)}\right)^k$ for all $k, n \in \mathbb{N}$, and we know that their limiting averages $\Lambda_k = \lim_{n \to \infty}(1/n)(n\rho_n)^k S_{n,k}$ exist and are finite by Assumption A2. Next, note that

$$\lim_{n \to \infty} \frac{n^s \rho_n^{2s}}{n} \times$$

$$\sum_{j=1}^{s} \sum_{\substack{1 \le i_1,\ldots,i_{s+1} \le n \\ \text{card}(\{i_1,\ldots,i_{s+1}\}) = j}} \left(\prod_{k=1}^{s+1} \left(w_{i_k}^{(n)}\right)^{d_k(T)}\right) = 0,$$

which is true since, for each $j \in [s]$, the term in the curly brackets above can be bounded as

$$0 \le \frac{n^{-s}(n\rho_n)^{2s}}{n} \sum_{\substack{1 \le i_1,\ldots,i_{s+1} \le n \\ \text{card}(\{i_1,\ldots,i_{s+1}\}) = j}} \left(\prod_{k=1}^{s+1} \hat{w}_n^{d_k(T)}\right)$$

$$\le \frac{1}{n^{s+1}} \binom{n}{j} (n\rho_n \hat{w}_n)^{2s} \le \frac{n^{j-s-1}}{j!} (n\rho_n \hat{w}_n)^{2s} = o(1),$$

where in the last equality we use the fact that Assumption A3 (iv): $n\rho_n \hat{w}_n = O(1)$. In other words, as $n \to \infty$ the summation over $j = 1, \ldots, s+1$ is dominated by its last terms $j = s+1$. We now take limits of both sides in (32) and use the preceding asymptotic dominance relation together with (33) to get that

$$\lim_{n \to \infty} \bar{m}_{2s}^{(n)} = \lim_{n \to \infty} (n\rho_n)^s \mu_{2s,s+1}$$

$$= \lim_{n \to \infty} \frac{(n\rho_n)^{2s}}{n^{s+1}} \sum_{T \in \mathcal{T}_{s+1}} \prod_{j=1}^{s} S_{n,j}^{r_j(T)}$$

$$= \lim_{n \to \infty} \sum_{T \in \mathcal{T}_{s+1}} \prod_{j=1}^{s} \frac{(n\rho_n)^{jr_j(T)}}{n^{r_j(T)}} S_{n,j}^{r_j(T)}$$

$$= \sum_{T \in \mathcal{T}_{s+1}} \prod_{j=1}^{s} \Lambda_j^{r_j(T)}, \tag{34}$$

where in the penultimate equality we have used the facts that, for any of the trees $T \in \mathcal{T}_{s+1}$, the degree distribution $(r_1, \ldots, r_s)$ satisfies $\sum_{d=1}^{s} d\, r_d = 2s$ and $\sum_{d=1}^{s} r_d = s+1$.

The last step of the proof relies on our ability to count the number of rooted ordered trees with any fixed degree distribution $\bar{r} = (r_1, \ldots, r_s)$, i.e. $\mathcal{T}_{s+1}(\bar{r})$. To determine $\text{card}(\mathcal{T}_{s+1}(\bar{r}))$, we rely on related properties of rooted ordered forests that are studied in [54, Chapter 5]. The total number of rooted ordered trees with degree distribution $\bar{r}$ is

given by:
$$\text{card}(\mathcal{T}_{s+1}(\bar{r})) = \frac{2}{s+1}\binom{s+1}{r_1,\ldots,r_s}. \quad (35)$$

We can now use (35) to transform the right-hand side of (34) into the expressions $m_{2s}$ claimed in the main theorem. Using the partition $\mathcal{T}_{s+1} = \bigcup_{\bar{r}\in\mathcal{R}_s}\mathcal{T}_{s+1}(\bar{r})$, we get

$$\lim_{n\to\infty}\mu_{2s,s+1} = \sum_{\bar{r}\in\mathcal{R}_s}\sum_{T\in\mathcal{T}_{s+1}(\bar{r})}\prod_{j=1}^{s}\Lambda_j^{r_j(T)}$$
$$= \sum_{\bar{r}\in\mathcal{R}_s}\text{card}(\mathcal{T}_{s+1}(\bar{r}))\prod_{j=1}^{s}\Lambda_j^{r_j}$$
$$= \sum_{\bar{r}\in\mathcal{R}_s}\frac{2}{s+1}\binom{s+1}{r_1,\ldots,r_s}\prod_{j=1}^{s}\Lambda_j^{r_j},$$

where (35) is invoked in the last equality, concluding the proof.
□

The results in Sections 7.1 and 7.2 enable us to claim that, under Assumptions A1 to A3, the spectral moments sequence $\{\mathbf{m}_k^{(n)}, k\in\mathbb{N}\}$ converges pointwise almost surely to the deterministic sequence $\{m_k, k\in\mathbb{N}\}$. Coup de grâce is to conclude the almost sure and weak convergence of the empirical spectral distributions from the pointwise almost sure convergence of their moments sequence; it is the province of the method of moments [64, Lemma B.3]: it is easy to verify that the limiting moment sequence $m_s, s\in\mathbb{N}$ satisfies Carleman's criterion, thus pointwise convergence of the moments sequence indeed implies convergence in distribution, completing the proof of the main result.

### ACKNOWLEDGMENTS

This work was supported in part by the United States National Science Foundation under grants CNS-1302222 and IIS-1447470.

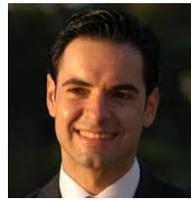

**Victor M. Preciado** received his Ph.D. degree in Electrical Engineering and Computer Science from the Massachusetts Institute of Technology in 2008. He is currently the Raj and Neera Singh Assistant Professor of Electrical and Systems Engineering at the University of Pennsylvania. He is a member of the Networked and Social Systems Engineering (NETS) program and the Warren Center for Network and Data Sciences. His research interests include network science, dynamic systems, control theory, and convex optimization with applications in socio-technical systems, technological infrastructure, and biological networks.

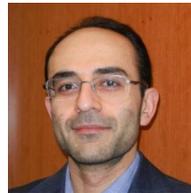

**M. Amin Rahimian** is a recipient of gold medal in 2004 Iran National Chemistry Olympiad. He was awarded an honorary admission to Sharif University of Technology, where he received his B.Sc. in Electrical Engineering-Control. In 2012 he received his M.A.Sc. from Concordia University in Montréal, and in 2016 he received his A.M. in Statistics from the Wharton School at the University of Pennsylvania, where he is currently a PhD student at the Department of Electrical and Systems Engineering and the GRASP Laboratory. He was a finalist in 2015 Facebook Fellowship Competition, as well as 2016 ACC Best Student Paper Competition. His research interests include network science, distributed control and decision theory, with applications to social and economic networks.